\documentclass[11pt,a4paper]{article}
\usepackage[fleqn]{amsmath}
\usepackage[cp1251]{inputenc}
\usepackage[english]{babel}
\usepackage{jcappub}
\usepackage{mathrsfs} 
\usepackage{hyperref}
\usepackage{xcolor}
\usepackage{cmap}
\title{Oleg Marichev: On the occasion of the eightieth anniversary}
\date{\today}

\author[1]{Askhabov S.N.}

\emailAdd{askhabov@yandex.ru}

\affiliation[1]{Chechen State University, Grozny, Russia} 
 
\author[2]{Kiryakova V.}

\emailAdd{virginia@diogenes.bg}

\affiliation[2]{Bulgarian Academy of Sciences, Sofia, Bulgaria}

\author[3]{Kravchenko V.V.}

\emailAdd{vkravchenko@math.cinvestav.edu.mx}

\affiliation[3]{Cinvestav, Campus Queretaro, Queretaro, Mexico}

\author[4]{Li C.P.}

\emailAdd{lcp@shu.edu.cn}

\affiliation[4]{Shanghai University, Shanghai, China}
	
\author[5]{Luchko Yu.F.}
	
\emailAdd{luchko@bht-berlin.de}
	
\affiliation[5]{Berlin University of Applied Sciences and Technology, Berlin, Germany}

\author[6]{Rogosin S.V.}

\emailAdd{rogosinsv@gmail.com}
		
\affiliation[6]{Belarusian State University, Minsk, Belarus}

\author[7]{Shishkina E.L.}

\emailAdd{elina.shishkina.math@gmail.com}
			
\affiliation[7]{Voronezh State University, Voronezh, Russia}
				
\author[8]{Sitnik S.M.}
				
\emailAdd{mathsms@yandex.ru}
				
\affiliation[8]{Belgorod State National Research University, Belgorod, Russia}
					
\author[9]{Vu Kim Tuan}

\emailAdd{vukimtuan1961@gmail.com}
					
\affiliation[9]{University of West Georgia, Georgia, USA}

\abstract{September 7, 2025 marked the 80th anniversary of the birth of Oleg Marichev. Marichev is known mathematician  which has developed many of Mathematica's algorithms for the calculation of definite and indefinite integrals and hypergeometric functions including Meijer $G$-function.}

\notoc

\makeatletter
\gdef\@fpheader{}
\makeatother

\begin{document}

\maketitle

\newpage
	\begin{figure}[h!]
		\centering
	\center{\includegraphics[width=0.5\linewidth]{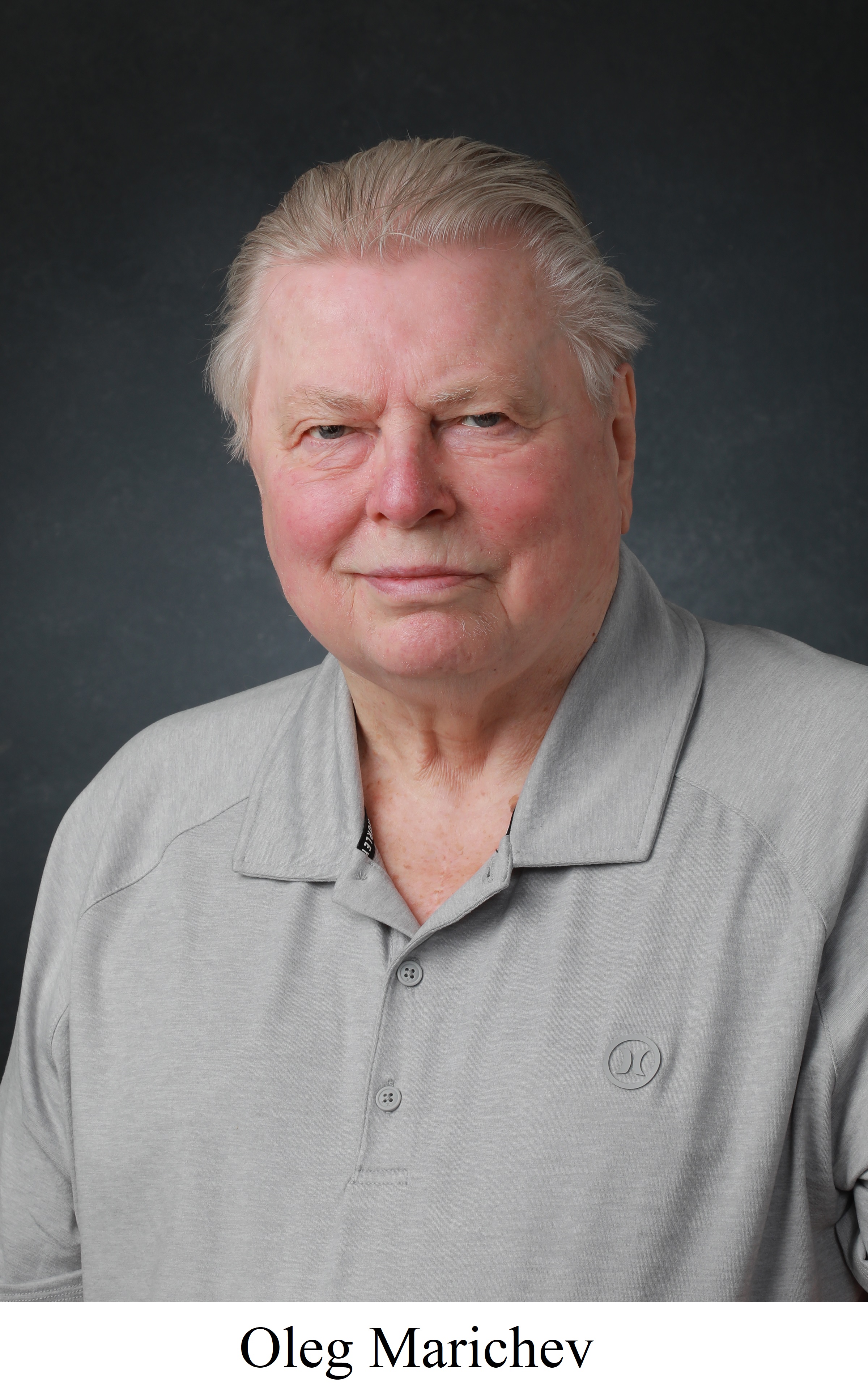} }
	%\caption{Oleg Marichev}
\end{figure}

Oleg Marichev was born on September 7, 1945, in the city of Velikiye Luki, located in the Pskov region of Russia.
In 1949, he and his family moved to Minsk, the capital of Belarus.
From 1952 to 1963, Oleg attended the 42nd Minsk school.
His interest in mathematics was sparked in the eighth grade by his teacher, Alexandra Bagreeva,
who introduced him to the method of mathematical induction and Newton's binomial theorem.
She also encouraged him to enroll in the school for young mathematicians at Belarusian State University (BSU).
It was there that he met his future wife, Anna, who was also a student at the school.
Oleg and Anna actively participated in mathematical olympiads, winning numerous prizes.
In 1963, Oleg graduated from school with a gold medal and was admitted to the Mathematics Department of BSU. After completing his university studies in 1968, he entered PhD course under the supervision of Professor Fyodor  Gakhov.
In 1966 Gakhov was elected a full member of the Academy of Sciences of the Byelorussian SSR.
Academician Gakhov proposed him to investigate the solvability of the mixed-type equations (in particular, Tricomi equation) which can be reduced to the singular integral equations.

While working on his PhD thesis, Oleg Marichev developed a keen interest in the field that later became the central focus of his scientific career, it is the theory of special functions. His fascination with generalized hypergeometric functions $\,_pF_q$, as well as Bessel, Legendre, Appell, Meijer, and Fox functions, was sparked by reading the book Generalized Hypergeometric Functions (1966) by Cambridge University professor Lucy Joan Slater. The application of integral transforms and special functions to solve differential equations formed the foundation of his PhD thesis, titled "Tricomi Boundary Value Problem for Some Equations of Mixed Type and Integral Equations with Special Functions in Kernels"\,, which he successfully defended in 1973. Following his dissertation defense, Oleg  Marichev was invited to join the faculty at Belarusian State University.
He became an outstanding mathematician who later gained recognition not only in the Soviet Union but also worldwide.

Oleg Marichev presented his initial studies on special functions and integral transforms, conducted between 1974 and 1978, in the monograph "Method of Calculating Integrals of Special Functions (Theory and Tables of Formulas)"\,. This book was published in Russian in Minsk in 1978 and later translated into English in 1983.

From 1968 to 1975 Oleg Marichev worked at the Department of Theory of Functions and Functional Analysis in the Faculty of Mathematics at Belarusian State University. One of his earliest students, Sergei  Rogosin  is now an associate professor in the Department of Analytical Economics and Econometrics at BSU. Another notable student of Oleg Marichev was I. E. Sadrigailo (1952--1976), who was awarded first prize in the All-Union Competition of Student Scientific Works.

In 1975, the Department of Theory of Functions and Functional Analysis was split into two separate departments: the Department of Theory of Functions and the Department of Functional Analysis. From 1975 to 1990, Oleg Marichev worked in the Department of Theory of Functions. He was awarded the title of associate professor in 1976 and was appointed to the position of acting professor.

The years 1981 to 1990 were a period of intensive scientific work for Oleg Marichev and his postgraduate students: Nguyen Thi Thanh, Vu Kim Tuan, Victor Adamchik, Semyon Yakubovich, Galina Grinkevich, and Nguyen Thanh Hai. Between 1984 and 1990, all of Oleg's PhD students successfully defended their PhD theses. Semyon Yakubovich defended his habilitation dissertation in 1996 in Minsk.  Notably, Vu Kim Tuan, after graduating BSU in 1984,  defended his PhD thesis in 1985 and habilitation dissertation in 1987 at the age of 26, having already authored more than 30 scientific articles in well-recognized journals by that time. Dr. Tuan later served as a distinguished chair professor at the University of West Georgia, USA, Dr. Victor Adamchik is a professor at the University of Southern California, USA, and Dr.  Semyon Yakubovich serves as a professor at the  University of Porto, Portugal. Associate Professor Galina Grinkevich spent her career at Vitebsk State University.

Since 1978, Oleg Marichev has been engaged in a large-scale project to compile integral tables more comprehensive than those Gradshteyn and Ryzhik, meticulously calculating and verifying thousands of complex integrals by hand. This effort culminated in the three-volume work Integrals and Series, co-authored with A. P. Prudnikov and Yu. A. Brychkov. In Russian, the three volumes are titled: Volume 1 --- Elementary Functions (1981), Volume 2 --- Special Functions (1983), and Volume 3 --- Additional Chapters (1986). Two more volumes (Vols 4 and 5) were published in English.

In 1988, a German translation of the first volume of Integrals and Series, supplemented with commentary, was published. Subsequently, in 1991, the first two volumes were published in Japanese in Tokyo. The algorithms developed in these books later served as the basis for creating functions within the Mathematica software system.

In 1990, Oleg Marichev defended his doctoral thesis titled Functions of Hypergeometric Type and Some of Their Applications to Integral and Differential Equations at Friedrich Schiller University of Jena in the former GDR. This university was the academic home of his longtime colleague and co-author, Professor H.-J. Glesske. Marichev's defense was among the first instances of a Soviet scientist defending a doctoral dissertation abroad. During this period, he also delivered a lecture at a conference in Leipzig, which caught the attention of a representative from Wolfram Research Inc. (WRI), an American company that he associates with since then.

While compiling tables of integrals and generalizing his research findings, Oleg Marichev aimed to automate algorithms for integral calculation. In 1980, he had the opportunity to pursue this goal: together with Ernst Krupnikov, he began implementing these ideas using the MIR single-user engineering computer in Novosibirsk. Together, they developed an integration program based on the application of the convolution theorem involving two Meijer $G$-functions. A few years later, during the perestroika period, the Department of Function Theory was able to establish a computer algebra laboratory.   At that time, Oleg Marichev, along with his former graduate student Viktor  Adamchik
and laboratory member Yuri Luchko (now full professor in mathematics at the Berlin University of Applied Sciences and Technology), created a prototype of the REDUCE automatic integration system within the computer algebra environment.

As a result, in 1990, Oleg  Marichev and Victor  Adamchik were invited to Champaign, USA, to demonstrate the capabilities of the REDUCE system at Wolfram Research. After showcasing their integral program, Stephen Wolfram offered to extend their visit by three months to experiment with implementing their algorithms in Mathematica. They continued working on the project and received temporary visas to stay for the following months. In March 1991, Oleg's wife Anna and their son Sergey, then eight years old, joined him in Champaign. In 1992, Oleg Marichev was granted permanent resident status.

From 1992 to 1997, Oleg Marichev worked actively at Wolfram Research Inc., concentrating on issues related to symbolic integration and the numerical analysis of the Meijer $G$-function, which is utilized in Mathematica under the name \texttt{MeijerG}. This function is recognized as one of the most complex and versatile special functions implemented in Mathematica.

Before developing a program to implement the Meijer $G$-function, Oleg  Marichev, at the special request of Stephen Wolfram, created programs for calculating indefinite integrals expressible in terms of elliptic integrals or other special functions, primarily those of hypergeometric type such as the Meijer $G$-functions.

Oleg Marichev contributed to the improvement and development of various Mathematica functions, including \texttt{Integrate} (in collaboration with Kelly Roach, Emily Martin, Alexey Bocharov, Victor Adamchik, and Daniel Lichtblau), \texttt{FunctionExpand} and \texttt{PowerExpand} (with Adam Strzebonski), \texttt{Series} and \texttt{Limit} (with Daniel Lichtblau), \texttt{BellY} (with Alexander Pavlik and Dan McDonald), \texttt{ContinuedFractionK} (with Charles Pugh and Michael Trott), among others. He received invaluable assistance in coding the \texttt{MeijerG} function from Jerry Kuiper (1953-995), a leading computational mathematics developer at Wolfram Research Inc., who tragically passed away in a car accident. Anna, Oleg's wife, also worked at the company from 1991 to 1997 and assisted him in testing integrals.

Oleg   Marichev has made  significant contributions to the development the theory of elementary and special mathematical functions, along with operations involving them such as various transformations, integration, and differentiation (including fractional)
in the Wolfram Mathematica system.
Many of Marichev's formulas exhibit conceptual differences from those published in existing literature; specifically, he accounted for the multivalued nature of functions of a complex variable and formalized them correctly. This approach enabled their accurate implementation in computer algebra systems, particularly within Mathematica.
For instance, the well-known classical formulas for the expansions of functions like $\ln(xy)$  and $(xy)^a$ were refined by incorporating adjustments based on
$\left[\frac{\pi-{\rm arg}(x)-{\rm arg}(y)}{2\pi}\right]$
where $[t]$ denotes the greatest integer not exceeding $t$. These modifications allowed for the numerical verification of all formulas presented in Mathematica. Such corrections arise from the fact that the arguments ${\rm arg}(x), {\rm arg}(y),...$, etc., of all variables are constrained within the range $-\pi{<}{\rm arg}(x){\leq}\pi$, $-\pi{<}{\rm arg}(y){\leq}\pi$, etc.
The implementation of functions such as \texttt{MeijerG}, the more general function \texttt{FoxH}, along with their numerous special cases; fractional calculus operators; and symbolic integration algorithms these are tools utilized daily by millions of users, ranging from schoolchildren to esteemed scientists.

Marichev transformed abstract theory into a practical tool, illustrating the systematic nature of the most complex calculated integrals. He demonstrated that, in the vast majority of cases, these integrals are merely special instances of broader "superintegrals"\, associated with "superfunctions"\,. The computation of these "superintegrals"\, can be achieved through the same "superfunctions"\, or their specific cases.

A significant role in the life of Oleg Marichev played  his work from 1984 to 1987 on the monograph "Integrals and Derivatives of Fractional Order and Some of Their Applications"\,, co-authored with Anatoly Kilbas (an Associate Professor in the Department of Function Theory at BSU, where Marichev was employed) and Professor Stefan Samko of Rostov State University.
Five years later, an expanded version of the book, featuring an additional 300 pages, was translated into English.
This handbook is well known as the Encyclopedia of the Fractional Calculus, contains the basic milestones of the theory of fractional order integrations and differentiation by  the 90s of 20th century, and brought to its 3 co-authors a worldwide recognition. Nowadays, in the era of booms of applications of this theory, almost no of the thousands such publications happen to miss to refer to it.
 In 2025, the book was translated into Chinese by Professor Changpin Li from the Department of Mathematics at Shanghai University and republished in China. This edition includes a substantial section titled "Overview of Fractional Calculus and its Computer Implementation in Wolfram Mathematica"\,, authored by Oleg Marichev and Elina Shishkina.

In 2018, following many years of collaborative work with Yu. A. Brychkov and N. V. Savishchenko, and incorporating a text authored by A. A. Kilbas based on the archival materials of A. P. Prudnikov, a reference book on the Mellin transform was published in English.

In the past 15 years, Oleg Marichev has been actively engaged in several projects. As of 2025, the Wolfram Functions website features approximately 307,000 formulas that describe the properties of functions, along with more than 31,000 formulas for 500 fundamental probability distributions (Wolfram Blog). Oleg Marichev provided solutions to four of the remaining 14 problems posed by Ramanujan, utilizing his programs for calculating continued fractions, specifically \texttt{ContinuedFractionK}, as discussed in the blog post "After 100 Years, Ramanujan Gap Filled"\,.
Most of the formulas available on the Wolfram Functions website are implemented in the Mathematica system using the Wolfram symbolic language, specifically through \texttt{MathematicalFunctionData} and \texttt{EntityValue} (for version 10.3, the implementation was carried out by Paco Jain and Michael Trott).

 The main viable result of Oleg Marichev is the deep developing of algorithm for evaluation integrals with its realization on computer. This algorithm covered near 80\% integrals, published in handbooks and articles. Oleg collected and developed largest at the world basic table of functions, combining pairs of them one can evaluate integrals from these pairs with their conditions of convergence and these pairs are the cases of more super generic functions with four lists of parameters (Meijer $G$-functions and Fox $H$-functions).

Congratulating Oleg Igorevich on his anniversary, we sincerely wish him good health, prosperity and new creative achievements!

\newpage
{\large{\bf References}}
	\begin{enumerate}
		%Гахов	
		\item		  {\bf Zverovich\,E.\,I., Komyak\,I.\,I.,  Mari\v{c}ev\,O.\,I.}  Fedor Dmitrievich Gahov (on his 70th birthday) (Russian).
		{\it Vestsi Akad. Nauk BSSR Ser. Fiz.-Mat.  Nauk}, 1976, vol.~1, pp.~131--133. \vspace{-3mm}
		
		\item  {\bf Kilbas\,A.\,A., Mari\v{c}ev\,O.\,I.}  Fedor Dmitrievich Gahov
		(Russian). {\it Vestsi Akad. Nauk BSSR Ser. Fiz.-Mat.  Nauk}, 1980, vol.~4, pp.~130--132. \vspace{-3mm}
		
		%вырождающиеся уравнения

		\item  {\bf Mari\v{c}ev\,O.\,I.} A mixed type equation with two lines of degeneracy in a nonsymmetrical domain (Russian).
		{\it Vestsi Akad. Nauk BSSR Ser. Fiz.-Mat.  Nauk}, 1969, vol.~6, pp.~74--78. \vspace{-3mm}

		\item  {\bf Mari\v{c}ev\,O.\,I.}   Boundary value problems for a mixed type equation with two curves of degeneracy (Russian).
		{\it Vestsi Akad. Nauk BSSR Ser. Fiz.-Mat.  Nauk}, 1970, vol.~5, pp.~21--29. \vspace{-3mm}
		
		\item  {\bf Mari\v{c}ev\,O.\,I.} Two Volterra equations with Horn functions. {\it Soviet Math. Dokl.}, 1972,
		vol.~13(3), pp.~703--707. \vspace{-3mm}

		\item  {\bf Mari\v{c}ev\,O.\,I.} Singular boundary value problems for a generalized biaxially symmetric  Helmholtz equation.
		{\it Soviet Math. Dokl.}, 1976,
		vol.~17(5), pp.~1342--1346. \vspace{-3mm}
		
		\item Scientific works of the Jubilee Seminar on Regional Problems, dedicated to the 75th anniversary of the birth of Academician of the Academy of Sciences of the BSSR F.D. Gakhov. Minsk~: Publishing House "University", 1985.\vspace{-3mm} %- 205 p.

		%Статьи 1974-1978 интегральные преобразования и специальные функции	
		
		\item  {\bf Mari\v{c}ev\,O.\,I.}  One class of integral equations of Mellin convolution type with  special functions in the kernels (Russian).
		{\it Vestsi Akad. Nauk BSSR Ser. Fiz.-Mat.  Nauk}, 1974, no.~1, pp.~126--127. Dep.in VINITI 21.09.73, no.7308, 18 pp. \vspace{-3mm}
		
		\item  {\bf Mari\v{c}ev\,O.\,I.} The Volterra equation of Mellin convolution type with $F_3$-function in the kernel (Russian).
		{\it Vestsi Akad. Nauk BSSR Ser. Fiz.-Mat.  Nauk}, 1974, no.~1, pp.~128--129. Dep. in VINITI 21.09.73, no. 7307, 23 pp. \vspace{-3mm}
		
		\item  {\bf Mari\v{c}ev\,O.\,I., Dang Din Kuang} The boundary value problem for the mixed type  equation of four order in the confounded domain with three  singular points.  (Russian).
		{\it Vestsi Akad. Nauk BSSR Ser. Fiz.-Mat.  Nauk}, 1974, no.~1, pp.~127--128. Dep. in VINITI no. 29, 35--75, 29 pp. \vspace{-3mm}
		
		\item  {\bf Mari\v{c}ev\,O.\,I.}  Some integral equations of Mellin convolution type containing the  special functions in the kernels (Russian).
		{\it Vestsi Akad. Nauk BSSR Ser. Fiz.-Mat.  Nauk}, 1976, no.~6, pp.~119--120. Dep.in VINITI 02.04.76, no.1640, 85 pp. \vspace{-3mm}

		\item  {\bf Mari\v{c}ev\,O.\,I.} Neumann and Dirichlet weight problems in the halfplane for the  generalized Euler-poisson-Darboux equation (Russian).  {\it Vestsi Akad. Nauk BSSR Ser. Fiz.-Mat.  Nauk}, 1976, no.~4, pp.~128--130. \vspace{-3mm}

		\item  {\bf Mari\v{c}ev\,O.\,I., Nguyen Ti Thanh}  Homogeneous Tricomi boundary value problem  for Lavrent'ev-Bitsadze equation in the strip similar to saw  with infinite  quantity of singular points (Russian).
		{\it Vestsi Akad. Nauk BSSR Ser. Fiz.-Mat.  Nauk}, 1976, no.~1, pp.~128--130. Dep. in VINITI, no. 29, 36--75, 48 pp. \vspace{-3mm}
		
		\item  {\bf Mari\v{c}ev\,O.\,I.} Full singular integral equations with power-logarithmic kernels  solved in the closed form (Russian).
		{\it  Vestnik Beloruss. Gos. Univ.}, 1978, vol. 1, no.~1, pp.~8--14. \vspace{-3mm}
		
		\item  {\bf Mari\v{c}ev\,O.\,I.} Integral operators with kernels generalizing the operators of  complex order (Russian).
		{\it Vestsi Akad. Nauk BSSR Ser. Fiz.-Mat.  Nauk}, 1978, no.~2, pp.~38--44.\vspace{-3mm}
		
		\item  {\bf Mari\v{c}ev\,O.\,I.} An integral representation of the solutions of a generalized biaxially symmetric Helmholtz equation and formulas for its inversion. {\it Differ. Uravn.}, 1978, vol. 14, no.~10, pp.~1824--1831.\vspace{-3mm}

		\item  {\bf Mari\v{c}ev\,O.\,I.}  Singular boundary value problems for Euler-Poisson-Darboux equation and  biaxially symmetric Helmholtz equation (Russian). {\it Proc. Vsesoyuz. conf. devoted  to 75th birthday of academician  I.G. Petrovskii}, MGU, Moscow, 1978, pp.~373--374.\vspace{-3mm}
		
		%желтая книга
		
		\item  {\bf Mari\v{c}ev\,O.\,I.} {\it Method of calculating integrals of special functions (theory and tables of formulas).} Minsk, Science and Technology, 1978.  \vspace{-3mm} %, 312 с.

		\item  {\bf Marichev\,O.\,I.} {\it Handbook of integral transforms of higher transcendental functions: theory and algorithmic tables.} Chichester, Ellis Horwood, 1983.\vspace{-3mm} % 336 p.
		
		\item  {\bf Prudnikov\,A.\,P., Brychkov\,Yu.\,A., Marichev\,O.\,I.} Calculation of integrals and Mellin transformation (Russian). {\it Mathematical analysis. Results of science and technology. USSR Academy of Sciences VINITI}, 1989, vol. 27, pp.~3--146.\vspace{-3mm}
		
		\item  {\bf Marichev\,O.\,I.}	 A method for calculating integrals of hypergeometric functions(Russian). {\it Dokl. Akad. Nauk BSSR}, vol.~25, no~7, pp. 590--593.\vspace{-3mm}
		
		\item  {\bf Marichev\,O.\,I.} Calculation of integral transformations of  hypergeometric functions (Russian). Generalized functions and their applications in mathematical physics. {\it Akad. Nauk SSSR. Vychisl. Tsentr. Moscow.}, 1981, pp. 323--331. \vspace{-3mm}
		
		%интегралы и ряды 1981-
		
		\item {\bf Prudnikov\,A.\,P., Brychkov\,Yu.\,A., Marichev\,O.\,I.}  {\it Integrals and series. Elementary functions (Russian).}  Moscow, Science, 1981.\vspace{-3mm} %, 800 с.
		
		\item {\bf Prudnikov\,A.\,P., Brychkov\,Yu.\,A., Marichev\,O.\,I.} {\it Integrals and series. Special functions (Russian).}  Moscow, Science, 1983. \vspace{-3mm} %, 752 с. 
		
		\item {\bf Prudnikov\,A.\,P., Brychkov\,Yu.\,A., Marichev\,O.\,I.}  {\it Integrals and series. Special functions. Additional chapters (Russian).}  Moscow, Science, 1986. \vspace{-3mm}
		
		\item {\bf Prudnikov\,A.\,P., Brychkov\,Yu.\,A., Marichev\,O.\,I.}  {\it Integrals and series. Elementary functions (Russian).} Second edition, Moscow, FizMatLit, 2002.\vspace{-3mm}
		
		\item {\bf Prudnikov\,A.\,P., Brychkov\,Yu.\,A., Marichev\,O.\,I.}  {\it Integrals and series. Special functions (Russian).} Second edition, Moscow, FizMatLit, 2003.\vspace{-3mm}
		
		\item {\bf Prudnikov\,A.\,P., Brychkov\,Yu.\,A., Marichev\,O.\,I.}  {\it Integrals and series. Special functions. Additional chapters (Russian).}  Second edition, Moscow, FizMatLit, 2003.\vspace{-3mm}

		\item  {\bf Prudnikov\,A.\,P., Brychkov\,Yu.\,A., Marichev\,O.\,I.} {\it Integrals and series. Vol. 1: elementary functions.}  New York, London, Paris, Montreux, Tokyo, Gordon and Breach, 1986.\vspace{-3mm}

		\item  {\bf Prudnikov\,A.\,P., Brychkov\,Yu.\,A., Marichev\,O.\,I.} {\it Integrals and series. Vol. 2: special functions.} New York, London, Paris, Montreux, Tokyo,  Gordon and Breach, 1986. \vspace{-3mm}

		\item  {\bf Prudnikov\,A.\,P., Brychkov\,Yu.\,A., Marichev\,O.\,I.} {\it  Integrals and series. Vol. 3: more special functions.} New York, NJ, Gordon and Breach, 1990. \vspace{-3mm}
		
		\item  {\bf Prudnikov\,A.\,P., Brychkov\,Yu.\,A., Marichev\,O.\,I.}  {\it Integrals and series. Vol. 4: direct Laplace transforms.} New York, NJ, Gordon and Breach, 1992. \vspace{-3mm}

		\item  {\bf Prudnikov\,A.\,P., Brychkov\,Yu.\,A., Marichev\,O.\,I.} {\it Integrals and series, Vol. 5: inverse Laplace transforms.} New York, NJ, Gordon and Breach, 1992. \vspace{-3mm}
		
		\item  {\bf K\"{o}lbig\,K.\,S.}  {\it Table errata: Integrals and series. Vol. 4.} [Gordon and Breach,   New York, 1992] by A. P. Prudnikov, Yu. A. Brychkov and O. I.  
		Marichev. Math. Comp. 1997. Vol.~66, No 220. P.~1765--1766.
		\vspace{-3mm}

		\item  {\bf Prudnikov\,A.\,P., Brychkov\,Yu.\,A., Marichev\,O.\,I.}  {\it Integrals and series. Vol.~1. Elementary functions. (Japaneese).}	Tokyo, Maruzen Co., Ltd, 1991. \vspace{-3mm}
		%	798 pp.
		
		\item  {\bf Prudnikov\,A.\,P., Brychkov\,Yu.\,A., Marichev\,O.\,I.}  {\it Integrals and series. Vol.~2. Special functions. (Japaneese).} 
		Tokyo, Maruzen Co., Ltd, 1991. \vspace{-3mm}
		% 752 pp.
		
		%неопределенные интегралы
		
		\item  {\bf  Brychkov\,Yu.\,A., Marichev\,O.\,I., Prudnikov\,A.\,P.,} 
		{\it Tables of indefinite integrals. (Russian).} Moscow, Nauka, 1986. \vspace{-3mm} %192 с. 
		
		\item  {\bf  Brychkov\,Yu.\,A., Marichev\,O.\,I., Prudnikov\,A.\,P.,}   
		{\it Tables of indefinite integrals. (Russian).} Moscow, FizMatLit, 2003. \vspace{-3mm} % 192 с. 
		
		\item  {\bf Brychkov\,Yu.\,A., Marichev\,O.\,I., Prudnikov\,A.\,P.} {\it Tables of indefinite integrals.}
		New York-London~: Gordon \& Breach Science Publishers/CRC Press, 1989. \vspace{-3mm}  %New York--London. ISBN 2 - 88124 - 710 - 5. (192 pages)
		
		\item  {\bf Brytschkow\,J.\,A., Maritschew\,O.\,I.,  Prudnikow\,A.\,P.} {\it Tabellen unbestimmter Integrale. (German).} Verlag, Harri Deutsch. 1992. \vspace{-3mm} % 192 pp.
		
		\item  {\bf Brychkov\,Yu., Marichev\,O., Savischenko\,N.}
		{\it Handbook of Mellin transforms.} New York, Chapman and Hall/CRC,  2018. \vspace{-3mm}
		
		%1982
		
		\item {\bf Marichev\,O.\,I.,  Tuan\,V.\,K.} Some properties of the $q$-gamma-function $\Gamma_q(z)$. {\it Dokl. Akad. Nauk BSSR}, 1982, vol.~26, no~6, pp. 488--491.\vspace{-3mm}
		
		\item {\bf Marichev, O. I.}   Conditions for the reversibility of Mellin-Barnes integrals.
		{\it Dokl. Akad. Nauk BSSR}, 1982, vol.~26, no~3, pp. 205--208.\vspace{-3mm}

		%1983 
		
		\item {\bf Marichev\,O.\,I.,  Tuan\,V.\,K.}  The problems of definitions and symbols of $G$- and $H$-functions of several variables. {\it Rev. T\'{e}cn. Fac. Ingr. Univ. Zulia. (Special Issue)} 1983, vol.~6. pp.~144--151.\vspace{-3mm}
		
		\item  {\bf Adamchik\,V.\,S.,  Marichev\,O.\,I.}  Representations of functions of hypergeometric type in  logarithmic cases (Russian).
		{\it Vestsi Akad. Navuk BSSR Ser. Fiz.-Mat. Navuk}, 1983, vol.~5, pp.~29--35. \vspace{-3mm}

		\item  {\bf Marichev\,O.\,I.}  Major terms of the asymptotic integral expansions of  the hypergeometric functions. {\it Short communications (Abstracts) of ICM-82. Warsazawa}, 1983, part~7, sect.~9, pp.~54. \vspace{-3mm}
		
		\item  {\bf  Marichev\,O.\,I.}  Asymptotic behavior of functions of hypergeometric type.
		{\it Vestsi Akad. Navuk BSSR Ser. Fiz.-Mat. Navuk}, 1983, vol.~4, pp.~18--25. \vspace{-3mm}

		\item {\bf  Marichev\,O.\,I., Vu Kim Tuan}  The definition of general $G$-function  of two variables, its principal cases and differential equations (Russian). {\it Differ. Uravn.}, 1983, vol.~19, no~10, pp. 1797--1799. Dep.in VINITI 08.02.83. no. 687--83, 20 pp.\vspace{-3mm}

		%1984
		
		\item {\bf Bry\v{c}kov\,J.\,A., Glaeske\,H.-J.,  Mari\v{c}ev\,O.\,I.}  Die Produktstruktur einer Klasse von Integraltransformationen. {\it Complex analysis and its  applications to partial differential equations.  Kongress und Tagungsberichte der Martin Luther Universitet  Halle Wittenberg}, 1984, pp.~59--60. \vspace{-3mm}
		
		\item {\bf Marichev\,O.\,I.}   On the representation of Meijer's $G$-function in the vicinity of singular unity. {\it  Complex analysis and   applications 1981, Varna.  Publ. House Bulgar.   Acad. Sci., Sofia}, 1984, pp.~383--398.  \vspace{-3mm}

		\item {\bf Marichev\,O.\,I.,  Kalla\,S.\,L.}  Behaviour of hypergeometric function $\,_pF_{p-1}(z)$
		in the vicinity of unity. {\it Rev. Tecn. Fac. Ingr. Univ. Zulia}, 1984, vol.~7, no~2, pp.~1--8. \vspace{-3mm}
		
		\item  {\bf Glaeske\,H.-J.,  Mari\v{c}ev\,O.\,I.} 
		The Laguerre transform of some elementary functions. {\it Z. Anal. Anwendungen}, 1984, vol.~3, no~3, pp.~237--244. \vspace{-3mm}
		
		\item  {\bf Adamchik\,V.\,S.,  Marichev\,O.\,I.}  Fundamental system of solutions for the operator $L_n=\left(z\frac{d}{dz}\right)^n+z$.
		{\it Differ. Uravn.}, 1984, vol.~20, no~6, pp. 1082--1083.\vspace{-3mm}

		\item {\bf Brychkov\,Y.\,A., Glaeske\,H.\,J., Marichev\,O.\,I.} Factorization of integral transformations of convolution type. {\it J. Math. Sci.}, 1985, vol.~30, pp. 2071--2094. \vspace{-3mm}

		%1985
		
		\item {\bf Marichev\,O.\,I.,  Tuan\,V.\,K.} Some Volterra equations with the Appell function $F_1$ 
		in the kernel (Russian).  {\it Proceedings of a commemorative seminar on  
			boundary value problems, Minsk}, 1985, pp. 169--172.\vspace{-3mm}

		\item {\bf Azamatova\,V.\,I.,  Marichev\,O.\,I.}
		The results of Belorussian mathematicians on integral equations   of convolution type and their applications  (Russian).  {\it Proceedings of a commemorative seminar on boundary value problems, Minsk}, 1985, pp. 27--30.  \vspace{-3mm}

		\item {\bf Marichev\,O.\,I.,  Tuan\,V.\,K.}  	 Composition structure of some integral transformations of   convolution type (Russian). {\it Reports of the extended sessions of  a seminar of the I. N. Vekua   Institute of Applied Mathematics. Tbilis. Gos. Univ., Tbilisi.},  1985, vol.~I, no. 1, pp. 139--142. \vspace{-3mm}
		
		%1986
		
		\item {\bf Marichev\,O.\,I.,  Tuan\,V.\,K.}  The factorization of G-transform in two spaces of functions. 
		{\it Complex analysis and applications. Varna~: Publ. House Bulgar. Acad. Sci., Sofia}, 1986, pp.~418--433. \vspace{-3mm}
		
		\item {\bf Brychkov\,Y.\,A., Marichev\,O.\,I.,  Yakubovich\,S.\,B.}  Integral Appell $F_3$-transformation with respect 
		to parameters. {\it Complex analysis and applications. Varna~: Publ. House Bulgar. Acad. Sci., Sofia}, 1986, pp.~135--140. \vspace{-3mm}
		
		\item  {\bf Brychkov\,Y.\,A., Marichev\,O.\,I., Yakubovich\,S.\,B.}  Factorization of integral transformations.
		{\it Publ. Math. Rep. Intern. Conference on Generalized Functions (Debrecen, 1984)}, 1986, vol.~33, pp.~161. \vspace{-3mm}
		
		\item {\bf Marichev\,O.\,I., Adamchik\,V.\,S., Tuan\,V.\,K.} 
		Solutions of a generalized hypergeometric differential equation. {\it Dokl. Akad. Nauk BSSR}, 1986, vol.~30, no~10, pp.~876--956. \vspace{-3mm}
		
		\item {\bf Bry\v{c}kov\,J.\,A., Glaeske\,H.-J.,  Mari\v{c}ev\,O.\,I.}  Die Produktstruktur einer Klasse von Integral transformationen. {\it Z. Anal. Anwendungen}, 1986, vol.~5, no~2, pp.~119--123. \vspace{-3mm}

		%1987 
		
		\item {\bf Marichev\,O.\,I.}  Compositions of fractional integrals and derivatives with power weights. Generalized functions, convergence structures and their applications. {\it Novi sad. Abstracts of conf. Dubrovnik}, 1987, pp.~37. \vspace{-3mm}
		
		\item  {\bf Yakubovich\,S.\,B., Tuan\,V.\,K., Marichev\,O.\,I.,  Kalla\,S.\,L.}
		A class of index integral transforms. {\it Rev. Tecn. Fac. Ingr. Univ. Zulia}, 1987, vol.~10, no~1, pp.~105--118. \vspace{-3mm}
		
		%1989
		
		\item {\bf Saigo\,M., Marichev\,O.\,I.,  H\`{a}i\,N.\,T.}  Asymptotic representations of Gaussian series $\,_2F_1$, Clausenian series 
		$\,_3F_2$, and Appell series $F_2$ and $F_3$ near boundaries of their convergence regions. {\it Fukuoka Univ. Sci. Rep.}, 1989, vol.~19, no~2, pp.~83--90.
		\vspace{-3mm}

		%1990
		
		\item {\bf Marichev\,O.\,I.}  Compositions of fractional integrals and derivatives  with power weights. {\it Fractional Calculus and its Applications. Proc. Inter. Conf. Nibon. Univ. (Tokyo, 1981), Ed. K. Nishimoto. Coll. Engin. Nihon Univ.}, 1990, pp.~94--99. \vspace{-3mm}
		
		%1991
		
		\item {\bf Vasiliev\,I.\,L.,  Marichev\,O.\,I.}  An operational method for the summation of series based on the   Laurent transform. {\it Dokl. Akad. Nauk BSSR}, 1991, vol.~35, no 4, pp.~296--299. \vspace{-3mm}
		
		%1992
		
		\item {\bf H\`{a}i\,N.\,T., Marichev\,O.\,I.,  Buschman\,R.\,G.}  Theory of the general H-function of two variables. {\it Rocky Mountain J. Math.}, 1992,  vol.~22, no~4, pp.~1317--1327. \vspace{-3mm}
		
		\item {\bf H\`{a}i\,N.\,T., Marichev\,O.\,I.,  Srivastava\,H.\,M.}  A note on the convergence of certain families of multiple   hypergeometric series. {\it J. Math. Anal. 
			Appl.}, 1992, vol.~164, no 1, pp.~104--115. \vspace{-3mm}

		%REDUCE system
		
		\item {\bf Adamchik\,V.\,S., Marichev\,O.\,I.}  The algorithm for calculating  integrals of hypergeometric type functions and its realization in REDUCE system. {\it Proc. Intern. Conf. ISSAC'90 (Tokyo,1990), New York, ACM}, 1990, pp.~212--224. \vspace{-3mm}
		
		%красная книга
		
		\item {\bf Samko\,S.\,G., Kilbas\,A.\,A., Marichev\,O.\,I.} {\it Integrals and derivatives of fractional order and some of their applications (Russian).} Minsk, Science and Technology, 1987.\vspace{-3mm} % 688 с.
		
		\item	{\bf Samko\,S.\,G.,  Kilbas\,A.\,A. and  Marichev\,O.\,L.}  {\it Fractional integrals and derivatives.} Amsterdam, Gordon
		and Breach Science Publishers, 1993.  \vspace{-3mm} 
		
		\item   {\bf Samko\,S.\,G., Kilbas\,A.\,A.,  Marichev\,O.\,I.}    {\it Fractional integrals and derivatives. Theory and Applications.} Translation on Chinese by prof. Changpin Li with " Overview of fractional calculus and its computer implementation in Wolfram Mathematica" by Marichev, O.I., Shishkina, E.L., China Science Publishing \& Media Ltd.(Science Press), Beijing, 1993.\vspace{-3mm}  %976 pp.
		
		\item {\bf Marichev\,O., Shishkina\,E.}
		Overview of fractional calculus and its computer implementation in Wolfram Mathematica. {\it Fractional Calculus and Applied Analysis}, 2024, vol.~27, pp.~1995--2062. \vspace{-3mm}

		%краевые задачи

		\item {\bf Marichev\,O.\,I., Kilbas\,A.\,A., Repin\,O.\,A.}  {\it Boundary value problems for partial differential equations with discontinuous coefficients (Russian).} Samara, Publishing house of Samara State Economic University, 2008. \vspace{-3mm}
		
		%Wolfram site
		
		\item http://functions.wolfram.com \vspace{-3mm}
		
		\item https://blog.wolfram.com/2008/05/06/two-hundred-thousand-new-formulas-on-the-web/ \vspace{-3mm}
		
		\item https://functions.wolfram.com/ElementaryFunctions/Log/16/04/01/0008/ \vspace{-3mm}
		
		\item {\bf Marichev\,O.\,I.}   Fractional Integro-Differentiation for Mathematica. {\it Book of Abstracts. International Symposium on Analytic Function Theory,  Fractional Calculus and Their Applications (In Honour of Professor H.M. Srivastava on his 65th Birth Anniversary, University of Victoria, PIMS}, 2005, pp.~33.\vspace{-3mm}

		%Теория вероятностей
		
		\item https://blog.wolfram.com/2013/02/01/the-ultimate-univariate-probability-distribution-explorer/ \vspace{-3mm}

		\item {\bf Brychkov\,Yu.\,A., Marichev\,O.\,I., Sofotasios\,P.\,H.}
		{\it Probability distribution functions and their characteristics (Russian)}, Handbook, Moscow, Max press, 2018. \vspace{-3mm}  %, 96 с.

		\item {\bf Kiryakova\,V., Kilbas\,A.\,A., Samko\,S.\,G.,  Marichev\,O.\,I.}  70th anniversary of Professor Rudolf Gorenflo. {\it  Fractional Calculus and Applied Analysis}, 2000,  vol.~3, no 4, pp.~339--342.\vspace{-3mm}

		\item https://functions.wolfram.com/Bessel-TypeFunctions/BesselJ/06/01/03/01/0006/ \vspace{-3mm} 
		
		\item https://functions.wolfram.com/ElementaryFunctions/Power/16/04/ \vspace{-3mm} 
		
		\item https://functions.wolfram.com/ElementaryFunctions/Sqrt/16/03/01/0007/ \vspace{-3mm} 
		
		\item https://functions.wolfram.com/ElementaryFunctions/Power/16/04/01/01/0008/ \vspace{-3mm} 
		
		\item {\bf Marichev\,O.\,I., Slavyanov\,S.\,Yu., Brychkov\,Yu.\,A.} Bell polynomials in the Mathematica system and asymptotic solutions of integral equations. {\it Theoret. and Math. Phys.}, 2019, vol.~201, no~3, pp. 1798–1807. \vspace{-3mm}
		
		\item  {\bf Yakovlev\,S.\,L., Andronov\,I.\,V., Suslina\,T.\,A.,
			Fedotov\,A.\,A.,  Its\,A.\,R.,  Motovilov\,A.\,K.,  Farafonov\,V.\,G.,  Kazakov\,A.\,Ya.}. In Memory of Sergei Yuryevich Slavyanov. 
		{\it Theor. Math. Phys.}, 2019, vol.~201, pp.~1543--1544.\vspace{-3mm}
		
		\item	 https://mathworld.wolfram.com/LauricellaFunctions.html \vspace{-3mm}
		
		\item https://blog.wolfram.com/author/michael-trott/ \vspace{-3mm}
		\item  {\bf Trott M.} {\it The Mathematica guidebook for programming.} 
		New York, Telos-Springer, 2004. \vspace{-3mm}
		\item  https://blog.wolfram.com/2013/05/17/making-formulas-for-everything-from-pi-to-the-pink-panther-to-sir-isaac-newton/ \vspace{-3mm}
		\item  https://blog.wolfram.com/2013/07/19/using-formulas-for-everything-from-a-complex-analysis-class-to-political-cartoons-to-music-album-covers/ \vspace{-1mm} 
	\end{enumerate}

% Main text
%\section{}\label{}

% Use the style of numbering in square brackets.
% If nothing is used, default style will be taken.
%\begin{enumerate}[a)]
%\item 
%\item 
%\item 
%\end{enumerate}  

%\begin{itemize}
%\item 
%\item 
%\item 
%\end{itemize}  

%\begin{description}
%\item[]
%\item[] 
%\item[] 
%\end{description}  

%\clearpage %%Remove this from your manuscript

% Figure
%\begin{figure}%[]
 % \centering
%    \includegraphics{}
 %   \caption{}\label{fig1}
%\end{figure}

%\begin{table}%[]
%\caption{}\label{tbl1}
%\begin{tabular*}{\tblwidth}{@{}LL@{}}
%\toprule
%  &  \\ % Table header row
%\midrule
% & \\
% & \\
% & \\
% & \\
%\bottomrule
%\end{tabular*}
%\end{table}

% Uncomment and use as the case may be
%\begin{theorem} 
%\end{theorem}

% Uncomment and use as the case may be
%\begin{lemma} 
%\end{lemma}


%\section{}\label{}


%\bibliographystyle{cas-model2-names}

%\bibliography{cas-refs}



\end{document}